\documentclass[12pt]{article} \newtheorem{Defn}{Definition} 
\newtheorem{lemma}{Lemma} 
\newtheorem{prop}{Proposition}\newtheorem{thm}{Theorem} 
\newcommand{\tensor}{\otimes} \title{Volume form as volume of 
infinitesimal simplices} \author{Anders Kock}

\begin{document}
\maketitle

\small
{\bf Abstract.}  In the context of Synthetic Differential Geometry, we 
describe the {\em square volume} of a ``second-infinitesimal simplex'', 
in terms of square-distance between its vertices.  The square-volume 
function thus described is symmetric in the vertices.  The 
square-volume gives rise to a characterization of the {\em volume 
form} in the top dimension.
\normalsize 

\section*{Introduction}

In the context of Synthetic Differential Geometry (SDG), it is possible to
express
some of the ``infinitesimal geometric'' notions for a manifold $M$
directly in terms of the points of the manifold, rather than in terms of the
tangent bundle $TM$; this was demonstrated for the case of connections,
differential forms, and
Riemannian metric, in a series of papers,  \cite{DFVG}, \cite{CTC},
\cite{CCBI}, \cite{Levi-Civita}, \cite{DFIC}, and the present note belongs to the same
series: here, we want to pass directly from the metric to the measure
of {\em volume}.

The crux is that in SDG, one has the notion of when points $x,y$ in a manifold
are first, resp.\ second, \ldots) order neighbours, (written $x\sim _1 
y$, resp.\ $x\sim _2 y$, \ldots ).

Recall \cite{CCBI} that an {\em infinitesimal $k$-simplex} is a $k+1$ 
tuple $(x_0 ,\ldots ,x_k )$ of mutual first order neighbours.  In this 
note, we shall consider also $k+1$ tuples $(x_0 ,\ldots ,x_k )$ of 
mutual {\em second} order neighbours, so we shall need a more elaborate 
terminology: call the former {\em first-infinitesimal} $k$-simplices, 
the latter {\em second-infinitesimal} $k$-simplices.

Recall also \cite{Levi-Civita} that a {\em Riemannian metric} 
synthetically may be expressed in terms of a function $g(x_0 ,x_1)$, 
applicable to any second-infinitesimal 1-simplex $(x_0 ,x_1 )$, and to 
be thought of as the {\em square length} (square of the length) of the 
simplex.  In particular, $g$ should be symmetric in its two arguments, 
and vanish on any first-infinitesimal 1-simplex.  Actually, the second 
condition follows from the first (provided $g(x,x)=0$), see 
Proposition \ref{one} below.

We shall in the present note construct, for any second-infinitesimal
$k$-simplex in a Riemannian manifold $(M,g)$, its {\em square volume}, and
prove, Proposition \ref{propsym}, that the square-volume is a symmetric
function in the $k+1$ vertices; also, the square volume will vanish if
two of the vertices are first order neighbours.  For the case where
$k=n$, the dimension $n$ of the manifold, it is possible in a certain
sense to extract square roots and derive a characterization of the
{\em volume form} (volume of first-infinitesimal $n$-simplices) of a 
Riemannian manifold; see Theorem \ref{volumeform}.

\section{Square volume: Heron's and Gram's constructions}

The set of first-infinitesimal 1-simplices in a manifold $M$ equals what in
other contexts is called {\em the first neighbourhood of the diagonal},
$M_{(1)}$, and similarly the set of second-infinitesimal simplices is the
second neighbourhood of the diagonal, $M_{(2)}$. As observed in \cite{DFVG}
Proposition 4.1, a function $f(x,y)$, vanishing on the diagonal, $f(x,x)=0$, 
is alternating on the first neighbourhood, $f(x,y)=-f(y,x)$.  For 
reference, we state this fact again, together with a slight extension 
of it:

\begin{prop} If a function in two variables vanishes on the diagonal, then it
is alternating ($f(x,y)=-f(y,x)$) on $M_{(1)}$.
If a function is symmetric ($f(x,y)=f(y,x)$) on
$M_{(2)}$, and vansihes on the diagonal, then it  vanishes on $M_{(1)}$.
	\label{one}\end{prop}
	For, on $M_{(1)}$, $f$ is alternating as well as symmetric.
	
	\medskip

For a Riemannian manifold, i.e.\ a manifold $M$ equipped with a 
Riemannian metric $g$, we shall construct ``square
$k$-volumes'' for any second-infinitesimal $k$-simplex (where $k \leq n$,
 $n$  the dimension of the manifold $M$ under consideration).

The basic idea for the construction of a square $k$-volume function
goes, for the case $k=2$, back to Heron of Alexandria (perhaps even
to Archimedes); they knew how to express the square of the area of a 
triangle $S$ in terms of an expression involving only the lengths 
$a,b,c$ of the three sides:
$$area^2 (S)  = t\cdot (t-a)\cdot (t-b) \cdot (t-c)$$
where $t= \frac{1}{2} (a+b+c)$. Substituting for $t$, and multiplying out,
one discovers (cf. \cite{Coxeter} \S 18.4) that all terms involving an odd
number of any of the variables
$a,b,c$ cancel, and we are left with
\begin{equation}area ^2 (S) = \frac{1}{16} (-a^4 -b^4 -c^4 + 2a^2 b^2 + 2 a^2
c^2 + 2 b^2
c^2 ),\label{zero} \end{equation}
 an expression that only involves the {\em
squares} $a^2$, $b^2$ and $c^2$ of the lengths of the sides.

For the case where the triangle is a second-infinitesimal 2-simplex,
the three first terms in the above parenthesis vanish, and we are
left with
$$area^2 (S) = \frac{1}{8} (a^2 b^2 +  a^2 c^2 +  b^2 c^2 ),$$
or expressed in terms of $g$ and the vertices $x_0 , x_1 , x_2 $ of
the simplex
\begin{equation} area ^2 ( x_0 , x_1 , x_2 ) = \frac{1}{8} (g(x_0 ,x_1 )
g(x_0 ,
x_2) +
g(x_1 ,x_0 ) g(x_1 , x_2 ) + g(x_2 , x_0 ) g(x_2 ,x_1 )).
\label{heron} \end{equation}
Since $g$ is symmetric in its two arguments, the expression here is 
symmetric in the three arguments.  Also, if two vertices, say $x_0$ 
and $x_1$ are equal, the two first terms in (\ref{heron}) vanish, but 
the third one does too, being $g(x_2 , x_0 )^2$.  (Note that the 
values of $g$ are always numbers with square zero.) From Proposition 
\ref{one} then follows that (\ref{heron}) vanishes if $x_0$ and $x_1$ 
are 1-neighbours; similarly, by symmetry, for $x_0 ,x_2$ or $x_1 , 
x_2$.

Since
the formula (\ref{heron}) only involves $g$, we do
have a ``square-area'' construction for second-infinitesimal
2-simplices, satisfying the conditions mentioned in the introduction:
symmetry, and vanishing if two vertices are 1-neighbours.

The expression in (\ref{heron}) has a rather evident generalization to
higher $k$, but we shall prefer to generalize a different (less symmetric) 
expression for 
square-area. This formulation involves the ``Gram determinant''. It
is well known from linear algebra books that the square of the volume of
the box
spanned by $k$ vectors $x_1 , \ldots ,x_k$ in ${\bf R}^n$ is the
Gram-determinant of the $n\times k$ matrix $X$ with the $x_i$'s as columns, 
where the {\em Gram determinant} of such matrix $X$ is the determinant of the  
$k\times k$ matrix
$X^{T} \cdot X$; or, equivalently, the matrix whose $ij$'th entry is $x_i
\cdot x_j$, where $\cdot$ denotes the standard inner product in
${\bf R}^n$.  Recall also that this inner product may be written in
terms of square-norm, by the polarization formula
$$x\cdot y = \frac{1}{2} (||x+y||^2 -||x||^2 - ||y||^2 ).$$

Thus the volume of the $k$-box spanned by the vectors $x_1 , \ldots
,x_k$ can be expressed in terms of square-norm. The volume of the
$k$-{\em simplex} spanned by some vectors (and 0) differs from the volume 
of the {\em box} spanned by the vectors by a factor $1/k!$.  This 
motivates the following construction of a square-$k$-volume of a 
second-infinitesimal $k$-simplex in a Riemannian manifold $M,g$:
\begin{equation}vol ^2 (x_0 , \ldots ,x_k ) = \frac{1}{(k!)^2} \mbox{det}
((x_i - x_0) \bullet
(x_j - x_0 )),\label{main}\end{equation}
where in turn
\begin{equation}
	(x_i - x_0) \bullet
(x_j - x_0 ) = \frac{1}{2} (-g(x_i , x_j ) + g(x_i , x_0 ) + g(x_j ,x_0
)),\label{polarization}\end{equation}
 Note that the minus signs
used on the left hand side are purely symbolical, and also a question
whether $\bullet$ is bilinear makes no sense.

We shall want to make sense to the expression (\ref{polarization}) 
without assuming $x_i \sim _2 x_j$, but assuming, as before, $x_i \sim 
_2 x_0$ and $x_j \sim _2 x_0$, whence at least $x_i \sim _4 x_j$.  
This depends on choosing an extension $\overline{g}$ of $g$ to the 
fourth neighbourhood of the diagonal $M_{(4)}$, so that the first term 
on the right hand side of (\ref{polarization}) makes sense.  The value of the 
$\bullet$ expression will in general depend on the choice of 
$\overline{g}$, but the dependence will be ``controlled'', in a sense 
that is best made precise by working in a coordinate neighbourhood 
around $x_0$:

The Riemannian metric $g$ on a neigbourhood of $0$ in $R^n$
 takes the form of a matrix product
  \begin{equation}g(x,y) = (y-x)^T 
\cdot G(x) \cdot (y-x)
\label{ddd}\end{equation}
 where $G(x)$ for each 
$x$ is a symmetric $n\times n$ matrix with positive determinant.

 Let $x\sim _2 0$ and $y\sim _2 0$. The value of $x -0 \bullet y-0$, as given by (\ref{polarization}) 
 with $x_0 =0, x_i = x $, and $x_j =y$ will depend on the extension 
 $\overline{g}$ chosen; however, since a different choice will mean a 
 difference for $\overline{g}(x,y)$ which is of order $\geq 3$ in 
 $y-x$, and since $x$ and $y$ are assumed to be $\sim _2 0$, this 
 difference will be of total order $\geq 3$ in $x, y$, and so may be 
 subsumed in the ``error term'' $\epsilon$ in the conclusion in the 
 following Lemma.  So in the proof of the Lemma, we may as well assume 
 that the expression for $\overline{g}$ is the expression (\ref{ddd}) 
 for all $x$ ,$y$.

\begin{lemma} Assume $x\sim _2 0$ and $y\sim _2 0$ (hence $x\sim _4 y$). Then
	$$x \bullet y = x \cdot G(0) \cdot y \; +\; \epsilon (x,y),$$
where $\epsilon (x,y)$ is of total degree $\geq 3$ in $x$ and $y$.  In 
particular, $x\bullet x = x^T\cdot G(0) \cdot x$.
	
\label{aa}\end{lemma}

{\bf Proof.} Let us calculate $2 x\bullet y$, i.e.  
$-\overline{g}(x,y)+g(0,x)+g(0,y)$.  We get
$$-(y-x)^T \cdot G(x) \cdot (y-x) + x^T\cdot G(0) \cdot x+ y^T\cdot 
G(0) \cdot y;$$
we Taylor expand $G(x)$ and thus rewrite the first term as $-(y-x)^T 
\cdot G(0) \cdot (y-x)$ plus a term $\epsilon $ which is of degree 
$\geq 1$ in $x$, and well as quadratic in $y-x$, and hence of total 
degree $\geq 3$ in $x,y$.  The three terms involving $G(0)$ now are, 
by the standard polarization formula, 2 times $x^T \cdot G(0) \cdot 
y$.  This proves the first assertion.

The assertion about $x\bullet x$ is then clear.

 \medskip
 
We consider now the determinant on the right hand side of (\ref{main}) 
in the coordinatized situation, i.e.\ with the metric given in 
coordinates by (\ref{ddd}).

\begin{prop}Consider a second-infinitesimal $k$-simplex $x_0 
,\ldots ,x_k$.  Then the  determinant $det((x_i - x_0 )\bullet (x_j -x_0 
))$ equals the determinant of the $k\times k$ matrix $X^T \cdot G(x_0 
) \cdot X$, where $X$ denotes the $n\times k$ matrix with the $x_i 
-x_0$ as columns.
\label{inv}\end{prop}

{\bf Proof.} The $ij$'th entry $(x_i - x_0 ) \bullet (x_j - x_0 )$ may 
by Lemma \ref{aa} be written
$$(x_i - x_0 )^T \cdot G(x_0 ) \cdot (x_j -x_0 ) + \epsilon 
_{i,j},$$
where $\epsilon _{i,j}$ is of total degree $\geq 3$ in $(x_i - x_0 ), 
(x_j - x_0 )$.  We claim that all the ``error terms'' $\epsilon 
_{i,j}$ get killed when expanding the determinant as a sum of products 
of the entries.  Note that the error terms only occur in off-diagonal 
entries, $i\neq j$.  We may write $\epsilon _{i,j}$ as 
$\epsilon ' _{i,j} + \epsilon ^* _{i,j}$, where $\epsilon ' _{i,j}$ 
is of degree $\geq 2$ in $x_i -x_0$ and $\epsilon ^* _{i,j}$ is of 
degree $\geq 2$ in $x_j -x_0$ (bidegrees $(3,0)$ or $(0,3)$ do not 
occur, since $x_i \sim _2 x_0$ for all $i$).  Now $\epsilon ' _{i,j}$, 
living in the $j$'th column, gets, in each of the terms of the 
deteminant expansion, multiplied by something from the $i$'th column.  
But each entry in the $i$'th column,  say
$$(x_l - x_0)\cdot G(x_0 ) \cdot (x_i - x_0 ) + \epsilon _{l,i}$$
is of degree $\geq 1$ in $x_i - x_0$ and therefore kills $\epsilon ' 
_{i,j}$.  Similarly,
$\epsilon ^* _{i,j}$ gets killed by being multiplied by anything in 
the $j$'th row. 

So the only factors left in the determinant expansion are the $(x_i - 
x_0)^T \cdot G(x_0 ) \cdot (x_j -x_0 )$, i.e.\ those of the matrix 
$X^T \cdot G(x_0 ) \cdot X$, as claimed.

\medskip

Whereas the ``Heron'' formulas (\ref{zero}) and (\ref{heron}) for square area
are evidently symmetric in the three vertices $x_0 ,x_1 ,x_2 $, the Gram
determinant formula (\ref{main}) gives a preferential status to the vertex
$x_0$. Nevertheless,

\begin{prop} The expression (\ref{main}) for square volume of a
second-infi\-ni\-te\-si\-mal $k$-simplex
is symmetric in the $k+1$ arguments $x_0 ,\ldots ,x_k$. Also, the expression
vanishes if two of the vertices $x_i$ and $x_j$ are 1-neighbours.
\label{propsym} \end{prop} 

{\bf Proof.} Symmetry in the arguments $x_1 ,\ldots
,x_k$ is clear. So it suffices to prove that 
\begin{equation}vol^2 (x_0
,x_1 , \ldots ,x_k ) = vol^2 (x_1 ,x_0 , \ldots ,x_k ).
	\label{symmetry}\end{equation}
By Proposition \ref{inv}, it suffices to prove that
\begin{equation}
det (X^T \cdot G(x_0 ) \cdot X) = det (Y^T \cdot G(x_1 ) \cdot 
Y)
\label{qqq}\end{equation} 
where $X$ (as above) is the $n\times k$ matrix with 
columns $x_i - x_0 $ ($i= 1,2,\ldots ,k$) and $Y$ similarly is the 
$n\times k$ matrix with columns $x_i - x_1$, ($i= 0, 2 ,\ldots ,k$).  
Now each of the $k!$ terms in the determinant $det (X^T \cdot G(x_0 ) 
\cdot X)$ is of degree $\geq 2$ in $x_1 -x_0$.  For, the entries of 
the first row of the matrix $X^T \cdot G(x_0 ) \cdot X$ are linear in 
$(x_1 - x_0)$, and likewise: each entry in the first column is linear 
in $(x_1 - x_0)$; and finally the entry in position $(1,1)$ is 
quadratic in $(x_1 - x_0)$.  Any of the $k!$ terms of the determinant 
expansion contains as well a factor from the first row, as one from 
the first column, and thus is of degree $\geq 2$ in $x_1 -x_0$.

Now since each term in the determinant expansion of $X^T \cdot G(x_0 ) 
\cdot X$ is of degree $\geq 2$ in $x_1 - x_0$, we may replace the 
$G(x_0 )$ by $G(x _1)$; for Taylor expanding $G(x_1 )$ from $x_0$ 
gives $G(x_0 )$ plus terms of degree $\geq 1$ in $x_1 -x_0$, but these 
terms get killed because our expressions are already of degree 
$\geq 2$ in $x_1 -x_0$.

Thus it suffices to prove that the $k\times k$ matrices $X^T \cdot G 
\cdot X$ and $Y^T \cdot G \cdot Y$ have the same determinant, where 
$G$ denotes $G(x_1 )$.  This is a matter of elementary linear algebra.  It is 
easy to see that one can obtain the matrix $Y$ from the matrix $X$ by 
elementary column operations: changing sign on a column, and 
adding a 
multiple of one column to another one.  This can be expressed matrix 
theoretically by saying $Y = X\cdot Z$ where $Z$ is a $k\times k$ 
matrix of determinant $\pm 1$.  Hence we have
$$det(Y^T \cdot G \cdot Y ) = det( (X\cdot Z)^T \cdot G \cdot (X \cdot Z))$$
 $$= det(Z^T \cdot (X^T \cdot G \cdot X)\cdot Z) = det(X^T \cdot G \cdot X),$$
 since $det (Z) = det (Z^T ) = \pm 1$, proving (\ref{qqq}), and hence 
 the symmetry assertion of the Proposition.

 To see that the expression (\ref{main}) vanishes if two vertices are
 1-neighbours, it suffices, by the symmetry already proved,
 and by Proposition \ref{one}, to prove that we get value zero if two
 vertices are equal.  This is clear from Proposition \ref{inv} , since 
 a determinant of form $X^T \cdot G \cdot X$ is zero if two columns of 
 $X$ are equal.
 
  \medskip

 {\bf Remark.} The above exposition does not, even in dimension 2, prove or
 argue that the Heron expression and the Gram expression agree; I suppose
 this is rather elaborate, or requires some more penetrating concepts.
 
 \medskip

\section{Volume form}

We are now going to compare square volumes with differential forms in the top dimension.
Differential forms  are defined
on first-infinitesimal simplices, and any square-volume function
applied to such simplex yields value 0. So no immediate comparison can be made.
The comparison therefore proceeds  via a
somewhat ad hoc notion of {\em extended $k$-simplex} and {\em extended 
$k$-form} $\Omega$.  

\begin{Defn} An {\em extended} $k$-simplex at $x_0$ is a 
$k+1$-tuple $(x_0 ,\ldots ,x_k ) $ of points with $x_i \sim _2 x_0$ 
for all $i= 1,\ldots ,k$.
\end{Defn}  
In particular, any 
second-infinitesimal simplex qualifies as an extended simplex.

Note that we do not assume $x_i \sim _2 x_j$ unless $i$ 
 or $j$ is 0, and so the notion of extended simplex is asymmetric.  
 (There are not $(k+1)!$ symmetries of an extended 
 $k$-simplex, but only $k!$.)

It is possible to define ``square volume'' not just for 
 second-infinitesimal $k$-simplices, but also for extended 
 $k$-simplices.   To define $vol^2 (X)$ of an extended 
 $k$-simplex $X$, we first choose an extension of $g$ to 
 $\overline{g}$, as above, so (\ref{polarization}) and hence 
 (\ref{main}) make sense, and serve as the definition of $vol^2 (X)$; 
 we have to argue that the value does not depend on the choice of 
 $\overline{g}$.  This follows because a different choice will, when 
 working in coordinates, give a difference in the $\epsilon$ part of 
 the $(x_i - x_0)\bullet (x_j - x_0 )$ only.  The same argument as used 
 in the proof of Propsition \ref{inv} then shows that all these $\epsilon$ 
 differences vanish when expanding the determinant.

\begin{Defn} An {\em 
extended $k$-form} $\Omega$ is a law which to any  extended 
$k$-simplex assigns a number; $\Omega $ should be alternating (change sign 
when swapping two entries $x_i$ and $x_j$ with $i$ and $j$ different, and 
$\geq 1$) and give value 0 if $x_i = x_j $ for some $i\neq j$.  
\end{Defn}
It makes 
sense to say that the extended $k$-form $\Omega$ is an extension of 
the $k$-form $\omega$.   --- The notion of extended form will only be 
used for the case $k=n$.  --- We have, for general $k$,

\begin{lemma}
Assume that the extended $k$-form $\Omega $ extends the non-degenerate
$k$-form $\omega$,
and similarly $\Omega '$ extends $\omega '$. If
$$\Omega (X) ^2 = \Omega ' (X) ^2$$
for every extended $k$-simplex $X$, then locally  $\omega = \omega '$
or $\omega = -\omega '$.\end{lemma}.

{\bf Proof.} We intend to prove $\omega (x_0 ,\ldots ,x_k )= \omega '(x_0
,\ldots
,x_k )$ for any $k$-simplex $(x_0 ,\ldots ,x_k )$. We introduce coordinates,
and assume $x_0 =0 \in R^n$. Then $\Omega$ as a function of the remaining $k$
variables $x_1 , \ldots ,x_k $ may be seen as a function
$D_2 (n) \times \ldots \times D_2 (n) \to R$ ($k$ copies of $D_2 (n)$). The
ring $B$ of such functions is a tensor product of $k$ copies of the ring
$A$ of
functions on $D_2 (n)$; this latter ring is graded,
$$A = A_0 \oplus A_1 \oplus A_2.$$ So $B$ is multi-graded; for simplicity
of notation,
let us do the case where $k=2$, so $B$ is bigraded,
$$B= \bigoplus _{i,j} A_i \tensor A_j .$$
The assumption that $\Omega$ vanishes if one of the $x_i$'s is 0  now implies
that $\Omega$ belongs to the summand
$\bigoplus _{i>0,j>0} A_i \tensor A_j$,
and so may be written
$$\Omega = \omega + H$$
where $\omega \in A_1 \tensor A_1$ and $H\in A_1 \tensor A_2 \oplus A_2
\tensor
A_1 \oplus A_2 \tensor A_2$.
Then for (bi-)degree reasons ($A_3 =0$), it follows that $\Omega ^2 =
\omega ^2$
in the bigraded ring $B$. The bilinear function which $\omega$ defines may be
identified with the given form $\omega$, and since this form is
non-degenerate, we may from $\Omega ^2 = \omega ^2$  conclude that 
$\omega = \pm \sqrt{\Omega ^2}$, so $\omega $ can, modulo sign, be 
reconstructed from $\Omega ^2$.  This gives the desired uniqueness 
assertion.

\medskip

A similar argument gives that if $\Omega $ and $\Omega '$ are 
two extended $k$-forms extending $\omega$, then $\Omega (X)^2 
= \Omega ' (X)^2$, for any extended $k$-simplex $X$.

\medskip

Let $Vol ^2$ (with capital $V$) denote the square-volume construction in the top
dimension $n$, given by the formula (\ref{main}); as noted 
above,   $vol^2 (X)$ makes invariantly sense even for an {\em 
extended} 
$k$-simplex $X$.

We shall prove

\begin{thm} Let $M, g$ be a Riemannian manifold of dimension 	$n$.  Then 
there exists locally, and uniquely up to sign, an $n$-form $\omega$, 
such that for some, or equivalently for any, extended $n$-form  $\Omega$ 
extending $\omega$,
  $$\Omega (x_0 , \ldots ,x_n )^2  = Vol ^2 (x_0 , \ldots 
  ,x_n).$$
\label{volumeform} \end{thm}

The $n$-form (unique up to sign) thus characterized, deserves 
the name {\em volume form} for the Riemannian manifold $M,g$.

\medskip

{\bf Proof.} From the Lemma  follows the local
uniqueness of volume forms, modulo sign. (So, if $M$ is connected, there
are at most two volume
forms $\omega$ and $-\omega$.) The {\em existence} of a volume form will
proceed in coordinates, and will ultimately not be very different from the
standard description in terms of square root of the determinant of the
symmetric
matrix representing $g$ in coordinates (as in the proof of Proposition 
\ref{propsym}).  The $n$ form $\omega$ we construct as a candidate for 
a volume form is then defined by 
\begin{equation}\omega (x_0 , \ldots ,x_n ) = 
\frac{1}{n!}\;  \sqrt{det G(x_0 )}\; det
(x_1 - x_0 , \ldots ,x_n - x_0 );\label{candidate}\end{equation}
 the same 
expression then defines a value for {\em any} $x_0 , \ldots ,x_n$, so 
in particular defines an extended $n$-form $\Omega$.  Then by linear 
algebra,
\begin{equation}\Omega (x_1 - x_0 , \ldots ,x_n - x_0 ) ^2 = \frac{1}{(n!)^2}
det(X^T \cdot G(x_0 ) \cdot X),\label{eee}\end{equation}
 where $X$ denotes the 
$n\times n$ matrix with columns the $x_i - x_0$.  

We now observe  that \begin{equation}det(X^T \cdot G(x_0 ) \cdot X) = det 
((x_i - x_0 )\bullet (x_j - x_0 )),\label{fff}\end{equation} which 
will then prove the existence of the volume form.  This is essentially 
immediate from Proposition \ref{inv}, because the proof about the 
cancellation of error terms depended on $x_i \sim _2 x_0$, but not on 
$x_i \sim _2 x_j$.

\end{document}